\documentclass[12pt,fleqn]{amsart}

\usepackage[all]{xy}
\usepackage{tikz-cd}
\usepackage{amscd}
\usepackage{setspace}
\usepackage{korpi-lib}

\usepackage{comment}

\usepackage{color}
\usepackage{enumerate}
\usepackage{amssymb}
\usepackage{hyperref}
\newcommand{\bmd}{d_{\rm BM}}
\newcommand{\sgn}{\protect{\rm sgn}}

\title[Banach-Mazur distance between some $C(K)$-spaces] {Bounds for Banach-Mazur distances\\ between some $C(K)$-spaces}
\author[M. Korpalski]{Maciej Korpalski}
\author[G.\ Plebanek]{Grzegorz Plebanek}
\address{Instytut Matematyczny, Uniwersytet Wroc\l awski, pl.\ Grunwaldzki 2, 50-384 Wroc\-\l aw\\ Poland}
\email{Maciej.Korpalski@math.uni.wroc.pl \\ Grzegorz.Plebanek@math.uni.wroc.pl}
\date{}

\thanks{The first author was partially supported by Młody Badacz grants funded by University of Wrocław}

\begin{document}

\begin{abstract}
We present several results providing lower bounds for the Banach-Mazur distance 
\[
\bmd\big( C(K), C(L)\big)
\]
between Banach spaces of continuous functions on compact spaces.
The main focus is on the case where $C(L)$ represents the classical Banach space $c$ of convergent sequences.
In particular, we obtain generalizations and refinements of recent results from \cite{GP24} and \cite{MP25}.

Currently, it seems that one of the most interesting questions is when $K = [0, \omega]$ is a convergent sequence with a limit and $L = [0,\omega]\times 3$ consists of three convergent sequences. In this case, we obtain
\[3.53125 \leq \bmd\big(C([0,\omega]\times 3),C[0,\omega]\big) \leq 3.87513\]
\end{abstract}
\maketitle

\section{Introduction}
Given two isomorphic Banach spaces $X$ and $Y$, their Banach-Mazur distance $\bmd(X,Y)$ is defined as the infimum of distortions $\|T\|\cdot\|T^{-1}\|$ taken over all isomorphisms $T:X\to Y$.
If we consider two compacta $K$ and $L$ and the Banach spaces $C(K),C(L)$ of real-valued continuous functions with the usual supremum norm, then $\bmd\big( C(K), C(L)\big)<2$ implies that $K$ and $L$ are homeomorphic and, consequently, $C(K)$ is isometric to $C(L)$. 
This was independently proved by Amir \cite{Am65} and Cambern \cite{Ca67}. The threshold $2$ is sharp --- see Cohen and Chu \cite{CC95} for a discussion of this phenomenon.

Recall that for compact metric spaces $K$ there is a sound isomorphic classification of $C(K)$, see Pe{\l}czy\'nski \cite{Pe68}. 
Below, we consider the space $[0,\omega]$, the simplest infinite compactum consisting of a converging sequence and its limit;
clearly, $C[0,\omega]$ represents the classical Banach space of converging sequences. 
Recall that $C(K)$ is isomorphic to $C[0,\omega]$ if and only if $K$ is a scattered compact space of height $<\omega$. 

There are few pairs of compacta $K$ and $L$ for which $\bmd\big( C(K), C(L)\big)$
is determined. 
A remarkable exception is provided by the following clean result.

\begin {theorem}\label{i:0}
The formula
\[\bmd\big(C([0,\omega]^m), C([0,\omega])\big)= m+\sqrt{(m-1)(m+3)},\]
holds for every $m\ge 1$.
\end{theorem}
Here, the upper bound was given by Candido and Galego \cite{CG13}*{Corollary 1.3 and Theorem 1.4(b)}, while Malec and Piasecki \cite{MP25} recently obtained the corresponding lower estimate.

Our main contribution to this topic is to present a relatively simple idea that yields lower bounds of
$\bmd\big( C(K), C(L)\big)$ for some pairs of $K,L$ --- this is described in Lemma \ref{main1}.
It is a refinement of methods of Gordon \cite{Go70}, Gergont and Piasecki \cite{GP24}, but as an application, we get substantially shorter argument leading to their results.
Actually, our proof gives the estimates from Theorem \ref{i:0} in a more general setting; see Theorems \ref{fa:1} and \ref{fa:2} for details.

The rest of the paper is devoted to the study of $\bmd\big( C(K), C[0,\omega]\big)$ where
$K=[0,\omega]\times k$ and $k$ is a natural number. We use the convention that $k=\{0,1,\ldots,k-1\}$
so that the space $K$ consists of $k$ copies of a convergent sequence. 
At first glance, this case may appear rather innocent. Moreover, Gordon \cite{Go70} proved that 
\[\bmd\big( C([0,\omega]\times 2), C[0,\omega]\big) = 3.\] 
However, when we pass from two to three copies of $[0,\omega]$, we face a problem that seems to be much harder: Gergont and Piasecki \cite{GP24} proved that
\[3.23 < \bmd\big( C([0,\omega]\times 3), C[0,\omega]\big) < 3.89\] 
by an involved argument and computer-aided calculations.
We show, using Lemma \ref{main1} again, that the lower bound here is at least $3.53$.
However, this requires solving several systems of linear inequalities, for which computer assistance again proves indispensable. 
The code used for this project is stored on GitHub at 
\[\text{\url{https://github.com/ememak/Bounds-for-Banach-Mazur-distance}.}\]
In Section \ref{k3} we also give a slightly improved upper bound of around $3.875$, which can possibly be optimal.

\section{Preliminaries}\label{p}

Throughout, $K$ and $L$ denote compact Hausdorff spaces.
Every Banach space of the form $C(K)$ is equipped with the supremum norm and
the dual space $C(K)^\ast$ is identified with the space $M(K)$ of all regular signed Borel measures on $K$ of finite variation. 
Given $f\in C(K)$ and $\mu\in M(K)$, we simply write $\mu(f)$ for $\int_K f \diff\mu$.

Suppose that $T:C(K)\to C(L)$ is an isomorphism such that $\|Tg\| \ge \|g\|$ for every $g\in C(K)$. 
For each $y\in L$ we denote by $\nu_y$ the signed measure on $K$ defined for $g\in C(K)$ by $\nu_y(g)=Tg(y)$; in other words, $\nu_y=T^\ast\delta_y$. 
In this setting, we note the following.

\begin{lemma}\label{p:1}
Measures $\nu_y$ for $y\in L$ form a 1-norming subset of $M(K)$. 
Moreover, for every $h\in C(L)$ there is $\vf\in C(K)$ such that $\nu_y(\vf)=h(y)$ for every $y\in L$.
\end{lemma}

\begin{proof}
If $g\in C(K)$ and $\|g\|=1$, then $\|Tg\|\ge\|g\|=1$, so there is $y\in L$ such that $|\nu_y(g)|=|Tg(y)|\ge 1$.
This means that $\{\nu_y:y\in L\}$ is a 1-norming set.

For any $h\in C(L)$, there is $\vf \in C(K)$ satisfying $T\vf=h$. 
Then $\nu_y(\vf)=T\vf(y)=h(y)$ for every $y\in L$, as required. 
\end{proof}

In the case of $L=[0, \omega]$, we write $\nu_i$ for $i\in [0,\omega]$ and denote $\nu_\omega$ simply by $\nu$. 
Note that $\nu_i\to\nu$ in the $weak^\ast$ topology of $M(K)$.

\begin{lemma}\label{p:1.5}
If $T:C(K)\to C[0,\omega]$ is a norm-increasing isomorphism, then $|\nu|(K)\ge 1$. 
\end{lemma}

\begin{proof}
By Lemma \ref{p:1}, there is $\vf\in C(K)$ such that $\nu_i(\vf)=1$ for every $i\le\omega$. Then $\|\vf\|\le 1$ and $|\nu|(K)\ge \nu(\vf)=1$.
\end{proof}

Later, it will be convenient to collect the following facts concerning the variation of a signed measure.

\begin{lemma}\label{p:2}
Let $K$ be any compact space and $\mu_i,\mu \in M(K)$.
\begin{enumerate}[(a)]
 \item If $\mu_i\to \mu$ in the $weak^\ast$ topology, then $|\mu |(C) \le \liminf_i |\mu_i|(C)$ for every clopen set $C\sub K$.
 \item If $\|\mu\|\le t$, $h$ is a norm-one measurable function and $B\sub A$ are two measurable sets, then $|\mu(h\chi_A)|\le c$ implies $|\mu(h\chi_B)|\le (t+c)/2$.
 \item Suppose that $h\in C(K)$, $\|h\| = 1$ and $h$ vanish outside a clopen set $C\sub K$. Further, let $z \in C$, $h(z)=1$. Writing $e$ for the characteristic function of $\{z\}$ we have
 \[ |\mu|(C)\ge 2|\mu(e)|-|\mu(h)|.\]
\end{enumerate}
\end{lemma}

\begin{proof}
For $(a)$ take any $\eps>0$ and a continuous function $g:C\to [-1,1]$ such that
$\mu(g)> |\mu| (C)-\eps$. Then 
\[ |\mu|(C) < \mu(g)+\eps = \lim_i \mu_i(g)+\eps \le \liminf_i |\mu_i| (|g|) +\eps\le \liminf_i |\mu_i | (C)+\eps.\] 
For clause $(b)$ note that
\[ -c\le \mu(h\chi_A)=\mu(h\chi_B)+\mu(h\chi_{A\sm B})\le c, \mbox{ and }\]
\[ -t\le \mu(h\chi_B)-\mu(h\chi_{A\sm B})\le t,\]
so $-t-c\le 2\mu(h\chi_B)\le c+t$, as required.

Part $(c)$ follows from
\[ 2|\mu(e)|-|\mu(h)|=|\mu(e)| +(|\mu(e)|-|\mu(h)|) \le |\mu(e)| +|\mu(e)-\mu(h)| \le \]
\[\le |\mu|(e)+|\mu|(|h-e|)=|\mu|(|h|)\le |\mu|(C).\]
\end{proof}

\section{Basic tool} \label{basic}
We consider here any scattered compact space $K$ of finite height, a norm-increasing isomorphism $T:C(K)\to C[0,\omega]$ and the associated measures $\nu_i, \nu \in M(K)$. We write $t=\|T\|$ for the norm of $T$.

In the proof below, as well as elsewhere, we use the asymptotic symbol $\lesssim$ in the following sense:
$a\lesssim b$ means that the real-valued functions $a$ and $b$ defined for $\eps>0$
satisfy
 $\lim_{\eps\to 0^+} a(\eps)\le \lim_{\eps\to 0^+} b(\eps)$.

\begin{lemma}\label{main1}
Suppose that $x\in K^{(1)}$, $C\sub K$ is a clopen set containing $x$ and $f\in C(K)$ satisfies
\begin{enumerate}[(i)]
 \item $t \ge s= \|Tf\|=\sup_i |\nu_i(f)|>1$;
 \item $f(x)=1=\|f\chi_C\|$.
\end{enumerate}
Then
\[\mbox{\rm (\ref{main1}.a)}\quad \limsup_{i\to\infty} |\nu_i|(C)\ge 2\frac{s-|\nu(f)|}{s-1}-|\nu(f\chi_C)|;\]
\[\mbox{\rm (\ref{main1}.b)}\quad t\ge 2\frac{s-|\nu(f)|}{s-1}-|\nu(f\chi_C)| +|\nu|(K\sm C).\]
\end{lemma}

\begin{proof} 
We first fix a sequence of isolated points $x_n\in C$ converging to $x$. 
Write $e_n\in C(K)$ for the characteristic function of $\{x_n\}$.

Fix $\eps>0$ and consider the functions
\[ g_n=\frac{1}{s+\eps} \cdot f +\left(1-\frac{1-\eps}{s+\eps}\right)\cdot e_n .\]
Note that $f(x_n)>1-\eps$ for large $n$ and then $\|g_n\|\ge g_n(x_n)\ge 1$.


As $\nu_i(f)\to \nu(f)$, there is $i_0$ such that for every $i\ge i_0$ we have $|\nu_i(f)-\nu(f)|< \eps$. 
Then fix $N$ such that for every $n\ge N$ and every $i<i_0$ we have $|\nu_i(e_n)|<\delta$,
where $\delta$ will be specified in a while.
We infer that for every $i<i_0$ and $n\ge N$,
\[ |\nu_i(g_n)| < \frac{s}{s+\eps}+\left(1-\frac{1-\eps}{s+\eps}\right)\cdot\delta < 1 \]
whenever $\delta$ is small enough.
In other words, we have checked that the initial measures cannot norm $g_n$ for large $n$. Hence we have checked that
\medskip

\begin{scclaim}
For every $n\ge N$ there is $i=i(n)\ge i_0$ such that $|\nu_i(g_n)|\ge \|g_n\|$.
\end{scclaim}

We can now perform the following approximate calculations:
\[ 1 \le |\nu_{i(n)}(g_n)| \lesssim \frac{|\nu(f)|}{s} + |\nu_{i(n)}(e_n) |(1-1/s), \mbox{ so}\]
\[ |\nu_{i(n)}(e_n)| \gtrsim \frac{s-|\nu(f)|}{s-1} (>0).\] 

Note that a measure of finite variation may have only finitely many big atoms; hence the sequence $i(n)$ is unbounded. 
Using Lemma \ref{p:2}(c) (with $h=f\chi_C$), we get
\[ |\nu_{i(n)} |(C)\gtrsim 2\frac{s-|\nu(f)|}{s-1}-|\nu(f\chi_C)|.\]
Since we started from an arbitrary $\eps>0$, the above asymptotic formula shows that $(a)$ holds. 
In turn, \ref{main1}(a) and Lemma \ref{p:2}(a) (applied for $K \sm C$) yield \ref{main1}(b) and the proof is complete. 
\end{proof}

\section{When \texorpdfstring{$K^{(2)}\neq\emptyset$}{K'' not empty}}\label{fa}

We apply Lemma \ref{main1} to estimate the Banach–Mazur distance between $C[0,\omega]$ and spaces of the form $C(K)$, where $K$ has nonempty higher derivatives.

\begin{theorem}\label{fa:1}
Let $K$ be a compact space such that $K^{(2)}\neq \emptyset$. Then
\[ \bmd\big( C(K), C[0,\omega] \big) \ge 2+\sqrt{5}.\]
\end{theorem}

\begin{proof}
Let $T:C(K)\to C[0,\omega]$ be a norm-increasing isomorphism and let $\nu_i$ be
corresponding measures on $K$ (see Section \ref{p}).
Writing $t=\|T\|$, we shall prove that $t=\sup_i |\nu_i|(K)\ge 2+\sqrt{5}$. 

Since $K$ is necessarily scattered, we may fix an isolated point $z$ of $K^{(2)}$ and a sequence $y_m\in K^{(1)} $ converging to $z$.
Write $\theta=\nu(\{z\})$; without loss of generality, we can assume that $\theta\ge 0$.
Choose a clopen set $A_0\sub K$ such that $A_0\cap K^{(2)}=\{z\}$
and $|\nu|(A_0)\approx \theta$.

\begin{scclaim}
There is a nonempty clopen set $A_1\sub A_0$ such that
\[\nu(A_1)\approx 0 \mbox{ and } |\nu_i(A_1)|\lesssim \frac{t+\theta}{2}\mbox { for every } i.\]
\end{scclaim}

To prove the claim, fix pairwise disjoint clopen sets $C_n\sub A_0$ such that $C_n\cap K^{(1)}= \{y_n\}$ whenever $y_n\in A_0$.

For any $\eps>0$ there is $i_0$ such that $|\nu_i(A_0)-\nu(A_0)|\le\eps$ for every $i\ge i_0$.
Then we can choose $A_1$ among the sets $C_n$ such that $|\nu_i(A_1)|<\eps$ for every $i<i_0$ and $|\nu|(A_1)<\eps$.
By Lemma \ref{p:2}, we have $|\nu_i(A_1)|\le (t+\theta+\eps)/2$ for every $i\ge i_0$
and this verifies the claim.
\medskip

We now apply Lemma \ref{main1} for $f=\chi_{A_1}$ with $s\approx (t+\theta)/2$: since
$\nu(A_1)\approx 0$ we get
\[ t\gtrsim 2\frac{(t+\theta)/2}{(t+\theta)/2-1} +|\nu|(K\sm A_1).\] 
Here $|\nu|(K\sm A_1)\gtrsim \max(1,\theta)$ so, finally,
\[ t\ge 2\frac{t+\theta}{t+\theta -2} +\max(1,\theta)= 2+ \frac{4}{t+\theta-2}+ \max(1,\theta).\] 
If $\theta\le 1$ then
\[ t\ge 2+ \frac{4}{t-1}+ 1 \mbox{ that is } t^2-4t-1\ge 0 \mbox{ implying } t\ge 2+\sqrt{5}.\]

Otherwise, we have $\theta>1$ and 
\[ t\ge 2\frac{t+\theta}{t+\theta -2} +\theta.\] 
{\sc WolframAlpha} says that $t\ge 2+ \sqrt{\theta^2+4}>2+\sqrt{5}$, and the proof is complete.
\end{proof}

Extending the argument above we prove the following general result which, in particular,
gives the lower bound needed for Theorem \ref{i:0}.

\begin{theorem}\label{fa:2}
Let $K$ be a compact space such that $K^{(m)}\neq \emptyset$ for $m \geq 2$. Then
\[ \bmd\big( C(K), C[0,\omega] \big) \ge m+\sqrt{(m-1)(m+3)}.\]
\end{theorem}

\begin{proof}
We follow here the notation from the beginning of the previous proof.
Take $z\in K^{(m)}$ and suppose, as before, that $\theta=\nu(\{z\})\ge 0$ and $A_0$ is chosen.
We first extend the claim from the previous proof.
\medskip

\begin{scclaim}
There are clopen sets $A_0\supseteq A_1\supseteq A_2\supseteq A_{m-1}$, where $A_{m-1}\cap K^{(1)}\neq\emptyset$ 
and indices $i_0 < i_1 < \ldots$ such that 
\begin{enumerate}[---]
 \item $|\nu|(A_1)\approx 0$ and $|\nu_i(A_1)|\approx 0$ for $i<i_0$;
 \item $|\nu_i(A_1)|\lesssim (t+\theta)/2$ for every $i_0\le i<i_1$;
 \item $|\nu_i(A_1)|\approx 0$ for every $ i\ge i_1$;
 \item $|\nu_i(A_2)| \lesssim (t-1)/2$ for every $ i_1\le i< i_2$;
 \item \ldots
\end{enumerate}
\end{scclaim}

Indeed, we choose $A_1$ as before but now we can assume that $z_n\in A_1\cap K^{(m-1)}$
so $z_1$ is a limit of points $y_n\in K^{(m-2)}$. 
We find that the corresponding sets $C_n\sub A_1$ and set $A_2$ to be one of them. 
The only difference is that now $|\nu_i|(A_1)\lesssim t-1$
for large $i$, so Lemma \ref{p:2} gives $|\nu_i(A_2)|\lesssim (t-1)/2$.
We can continue in this manner until we reach the first derivative of $K$.
\medskip

We then consider a function which is a convex combination of the form
\[ f=p \cdot \chi_{A_1}+\frac{1-p}{m-2}\cdot \sum_{j=2}^{m-1}\chi_{A_j}.\]
We have $\|f\|=1$ and we want to choose $p\in (0,1)$ to minimize $\|Tf\|$.
This $p$ is determined by the equation
\[ p\cdot \frac{t+\theta}{2}=\frac{1-p}{m-2}\cdot \frac{t-1}{2}.\]
Then
\[ p=\frac{t-1}{(t+\theta)(m-2)+t-1}, \mbox{ so}\]
\[ \|Tf\|\lesssim s(\theta):=\frac{t-1}{(t+\theta)(m-2)+t-1}\cdot \frac{t+\theta}{2}.\]

Consider first the case $\theta\le 1$: Lemma \ref{main1} gives
\[ t\ge \frac{2s(\theta)}{s(\theta)-1} +1,
\]
and, for $m \geq 2, \theta \geq 0$, this implies\footnote{using {\sc WolframAlpha}}
\[t\ge(1/2)\Big(\sqrt{4m^2+4m(\theta+1)+\theta^2-6\theta-7}+2m-\theta+1\Big):=h(\theta).\]

We claim that the function $h(\theta)$ on the right hand side is decreasing on $[0,1]$. Indeed, the derivative
\[h'(\theta) = \frac{1}{2} \left(-1 + \frac{(-3 + 2 m + \theta)}{\sqrt{-7 + 4 m^2 - 6 \theta + \theta^2 + 4 m (1 +\theta)}}\right),\]
is negative --- if we suppose that $h'(\theta) \geq 0$, then
\[(-3 + 2 m + \theta)^2 \geq -7 + 4 m^2 - 6\theta + \theta^2 + 4 m (1 + \theta)\]
so $16 - 16 m \geq 0$,
which is a contradiction. We thus conclude
\[t\ge h(1)=(1/2)\Big(\sqrt{4m^2+8m-12}+2m\Big)=m+\sqrt{(m-1)(m+3)}.\]

Suppose now that $\theta\ge 1$; then
\[ t\ge \frac{2s(\theta)}{s(\theta)-1} +\theta.
\]
Analysis using {\sc Mathematica} (see the file "Derivatives.nb" in the \href{https://github.com/ememak/Bounds-for-Banach-Mazur-distance}{project files}) shows that $s'(\theta)\ge 0$ for $t\geq 3$. 
Hence, the function $s(\theta)$ is increasing, so
\[ t\ge \frac{2s(1)}{s(1)-1} +1,\]
giving the same estimate.
\end{proof}

The previous result may be slightly generalized. Consider a space $L$ such that $L^{(2)}=\emptyset\neq L^{(1)}$.
Then we can identify $L$ with $[0,\omega]\times \{0,1,\ldots, k-1\}$ for some $k$.
Then one can re-examine the proof of \ref{fa:2}: we have to deal with a finite number of sequences of converging measures
$\nu_i\to \nu$ but the essence will be the same.

\begin{corollary}\label{fa:4}
Let $K$ be a compact space such that $K^{(m)}\neq \emptyset$. Then
\[ \bmd\big( C(K), C(L) \big) \ge m+\sqrt{(m-1)(m+3)}\]
for every compact space $L$ with an empty second derivative.
\end{corollary}

It also seems worth noting what Theorem \ref{fa:1} and Corollary \ref{fa:4} mean for $K$ of infinite height. 

\begin{corollary}
If $K^{(\omega)} \neq \emptyset$ and $L^{(2)} = \emptyset$, then $\bmd\big(C(K), C(L)\big) = \infty$.
\end{corollary}

\section{When \texorpdfstring{$K^{(2)}=\emptyset$}{K'' is empty}}\label{k'' nonempty}
We compare here two compacta
\[ K=[0,\omega] \times k \mbox{ and } L=[0,\omega],\]
where $k$ (=$\{0,1,\ldots k-1\}$) is a fixed natural number $k\ge 2$. 
We again fix a norm-increasing isomorphism $T: C(K)\to C(L)$ and set $t=\|T\|$. 
We still write $\nu_i=T^\ast\delta_i$ for every $i\in L$ and denote $\nu_\omega$ simply by $\nu$.

We now decompose the limit measure $\nu$ as
\[ {\rm (D)}\qquad \nu=\sum_{m<k}\theta_m\delta_{(\omega,m)}+\nu',\]
where $\nu'$ vanishes at all endpoints $(\omega,m)$. Recall that $\nu_i\to\nu$ in the
$weak^\ast$ topology of $C(K)$. 
Lemma \ref{p:1.5} implies that
\[ |\nu'|(K)+\sum_{m<k} |\theta_m|\ge 1.\]

Let us first explain why we can in fact assume that $\theta_m\ge 0 $ for every $m<k$.

\begin{lemma} \label{main:2 Composition with an involution}
Let $T : C(K) \to C(L)$ be a norm-increasing isomorphism and $\sigma \in C(K)$ be a function such that $\sigma^2 = 1$.
Then $\wh{T}$ defined as $\wh{T}(g) = T (g \cdot \sigma)$ for $g\in C(K)$ is a norm-increasing isomorphism of the same norm.
\end{lemma}

\begin{proof}
Clearly, for every $f\in C(K)$ we have $\| \sigma\cdot f\|=\|f\|$ so 
\[ \|f\|=\| \sigma\cdot f \| \le \|T(\sigma\cdot f)\|=\| \wh{T}f\|, \mbox{ and}\]
\[ \|\wh{T}f\| = \| T(\sigma\cdot f)\|\le \|T\|\|\sigma\cdot f\|=\|T\|\|f\|.\]

The operator $\wh{T}$ is surjective: for any $h\in C(L)$ there is $f\in C(K)$ such that $Tf=h$; then $\wh{T}(\sigma\cdot f)=T(\sigma\cdot\sigma\cdot f)=Tf=h$.
\end{proof}\begin{corollary}\label{s:2}
Let $T: C(K)\to C(L)$ be a norm-increasing isomorphism and let $\sigma\in C(K)$ be a function such that $\sigma(n,m)={\rm sgn}(\theta_m)$ for every $n\le\omega$ and $m<k$.

Then $\wh{T}:C(K)\to C(L)$ defined as $\wh{T}f=T(\sigma\cdot f)$ for $f\in C(K)$ is also a norm-increasing isomorphism with the same norm.
Moreover, if we decompose the measure $\wh{\nu}=\wh{T}^\ast\delta_\omega$ as in (D) then
$\wh{\theta}_m\ge 0$ for every $m<k$. 
\end{corollary}

\begin{proof}
Since $\sigma \cdot \sigma = 1$, by Lemma \ref{main:2 Composition with an involution} it is enough to note that $\wh{\nu}$ restricted to a given level $(\omega+1)\times \{m\}$ is equal to ${\rm sgn}(\theta_m)\nu$. 
\end{proof}

It will be convenient to use the following notation. For any $n\in\omega$ and $m<k$ we write
\[A_m(n)=[n,\omega]\times \{m\}.\] 
 
\begin{corollary}\label{sc:1}
If $\theta_m=0$ for some $m<k$, then $t\ge 2+\sqrt{3}$.
\end{corollary}

\begin{proof}
Suppose, for instance, that $\theta_0=0$. Then $|\nu'|(A_0(n)) \approx 0$ for $n$ large enough, hence $\nu(A_0(n))\approx 0$. 
We apply Lemma \ref{main1} with $C=A_0(n)$ and $s=t$: (\ref{main1}.b) gives $t\gtrsim 2t/(t-1)+1$, so $t\ge 2+\sqrt{3}$.
\end{proof}


\begin{corollary}\label{sc:2}
If $k\ge 2$ and $I\sub k$ is any doubleton then
\[ t\ge \frac{2t}{t-1}+ |\nu'|(K)+ \sum_{m\notin I}\theta_m.\]
\end{corollary}

\begin{proof}
Suppose, for instance, that $I=\{0,1\}$ and $\theta_0\le\theta_1$. Consider the function
\[ f=\chi_{A_0(n)}- \big({\theta_0}/{\theta_1}\big) \chi_{A_1(n)},\]
where $n$ is large enough. We have $\nu(f)\approx 0$, so Lemma \ref{main1}, applied with $C=A_0(n)\cup A_1(n)$, gives the declared formula.
\end{proof}
The following was proved by Gergont and Piasecki \cite{GP24}*{Theorem 3.2}: 

\begin{theorem}\label{sc:3}
Given any $k\ge 2$,
\[ \bmd\big(C([0,\omega]\times k),C[0,\omega]\big) \ge \frac{\sqrt{3k^2-2k+1}+2k-1}{k}.\] 
\end{theorem}

\begin{proof}
Write $a=|\nu'|(K)$ and $b=\sum_{m<k} \theta_m$ for simplicity.
Let $c(I)=\sum_{m\in I} \theta_m $ for any $I\in[k]^2$. Note that
\[ \sum_{I\in [k]^2} c(I)=(k-1)b,\]
so there is $I$ such that 
\[ c(I)\le \frac{k-1}{\binom{k}{2}} b= \frac{2b}{k}.\]

We apply Corollary \ref{sc:2} for such $I$: we have $a+b\ge 1$ and hence 
\[ t\ge \frac{2t}{t-1}+a+\frac{k-2}{k}b\ge \frac{2t}{t-1}+\frac{k-2}{k},\]
and the assertion follows by solving the related inequality.
\end{proof}

Theorem \ref{sc:3} says, in particular, that
\[\bmd\big( C([0,\omega]\times 2), C[0,\omega]\big) \ge 3,\] 
which is the optimal bound, see Gordon \cite{Go70}.

The next lemma will be crucial for the next section. 

\begin{lemma}\label{sc:4}
Denote $d_m=\chi_{\{(\omega,m)\}}$ for $m<k$ and let $c$ be a constant such that $t/2\le c\le t$. Further let $J$ be a set of those
$m<k$ for which there is $i(m)$ such that $| \nu_{i(m)}(d_m)| > c$.
Then
\[ |\nu'|(K)+\frac{t-c-1}{c}\sum_{m\in J} \theta_m+\sum_{m\notin J} \theta_m \ge 1.\]
\end{lemma}

\begin{proof}
Note that the mapping $J\ni m\mapsto i(m)$ is injective since $c\ge t/2$.
We pick a function $\vf\in C(K)$ such that 
\[ \nu_{i(m)}(\vf)= -\sgn(\nu_{i(m)}(d_m)) \mbox{ for }m\in J,\]
while $\nu_i (\vf)=1$ for $i\notin\{i(m): m\in J\}$.
Note that then $\|\vf\|\le 1$; write $\vf_m= \vf(\omega,m)$ for simplicity. 

Suppose for a while that $i(0)=0$ and $\vf_0\ge 0$. 
If $\nu_0(d_0)>c$ then 
\[ -1=\nu_0(\vf)=\nu_0(d_0)\vf_0+ \nu_0(\vf- d_0)\ge \vf_0 c-(t-c),
\mbox{ so } \vf_0\le \frac{t-c-1}{c}.\]
If $\nu_0(d_0)<-c$ then 
\[ 1=\nu_0(\vf)=\nu_0(d_0)\vf_0+ \nu_0(\vf- d_0)\le \vf_0 (-c)+t-c,\]
so, again, $\vf_0\le (t-c-1)/c$.

Recall that $\theta_m\ge 0$ for every $m<k$. Applying the remark above and writing $I=\{m<k: \vf_m\ge 0\}$, we conclude that
\[ 1=\nu(\vf)=\nu'(\vf)+\sum_{m<k} \vf_m\theta_m\le |\nu'|(K)+\sum_{m\in I} \vf_m\theta_m \le\] 
\[\le |\nu'|(K)+\sum_{m\in I\cap J} \vf_m\theta_m+\sum_{m\in I\sm J}\vf_m\theta_m \le\]
\[ \le |\nu'|(K)+ \frac{t-c-1}{c} \sum_{m\in I\cap J} \theta_m+\sum_{m\in I\sm J} \theta_m
\le |\nu'|(K)+ \frac{t-c-1}{c} \sum_{m\in J} \theta_m+\sum_{m\notin J} \theta_m ,\]
and we are done.
\end{proof}


\section{The mysterious case \texorpdfstring{$k=3$}{k = 3}}\label{k3}


The result of Gergont and Piasecki reproduced here as Theorem \ref{sc:3} states that

\[ \bmd\big(C([0,\omega]\times 3),C[0,\omega]\big) \ge \frac{\sqrt{3\cdot 3^2-2\cdot 3+1}+2\cdot 3-1}{k}=\frac{\sqrt{22}+5}{3}\approx 3.23.\] 

We outline here the method of proving that the distance in question is actually greater than 3.5 and enclose an analysis of a result from \cite{GP24} giving its upper bound.

\subsection{Lower bound}
Following the notation of the previous section, we additionally assume that $0\le \theta_0\le \theta_1\le\theta_2$.
Write $a=|\nu'|(K)$, $b=\theta_0+\theta_1+\theta_2$; note that $\theta_2\ge b/3$.

The main idea is to apply Lemma \ref{main1} for a number of functions $f\in C(K)$ and use Lemma \ref{sc:4} for
a certain constant $c\ge t/2$ to formulate a system of linear inequalities in the variables $\theta_0,\theta_1,\theta_2,a\ge 0$ with $t$ as a parameter. Lemma \ref{sc:4} defines a set $J\sub \{0,1,2\}$ --- depending on its structure we get four different cases.
We ask the following question: \textit{what is the maximal value of $t$ for which none of these systems of inequalities has a solution?}
This argument will produce a certain lower bound for $t$. 

Recall that a system of linear inequalities in four variables defines a polytope in ${\mathbb R}^4$ and it is nonempty if and only if it has a vertex.
Such a vertex is uniquely determined by four linearly independent equations related to those inequalities.
Hence, a manual analysis is, in principle, possible; however, in our situation, the computations are too involved, and we had to rely on computational assistance.

 Let us briefly explain the origin of these inequalities.
For every $n$ and $m<3$ we write $A_m(n)=[n,\omega]\times\{m\}$.
Note that, given $\eps>0$, there is $n_0$ such that for every $m<3$ we have
$|\nu'|(A_m(n_0))<\eps$. Then there is $i_0$ such that for every $i\ge i_0$
\[ |\nu_i( A_m(n_0))-\theta_m|<\eps.\]

Consider the function
\[ f=\chi_{A_0(n_0)}- \big({\theta_0}/{\theta_1}\big) \chi_{A_1(n_0)},\]
where $n_0$ is large enough --- a norm-one function for which $\nu'(f)\approx 0$ and $\nu(f)\approx 0$. 
Hence, Lemma \ref{main1} applied for $C=A_0(n_0)\cup A_1(n_0)$ and $s=t$ gives 
\[ t\gtrsim \frac{2t}{t-1}+\theta_2+a.\] 
This will give a suitable bound if $\theta_2+a$ is large enough.

Then consider
\[ f=\chi_{A_0(n_0)},\]
which will give the estimate
\[ t\ge 2\frac{t-\theta_0}{t-1}-\theta_0+\theta_1+\theta_2 +a,\] 
a suitable one whenever $\theta_0$ is small.


Other inequalities refer to the set $J\sub \{0,1,2\}$ defined in \ref{sc:4}; 
suppose, for instance, that $0\notin J$.
 
We choose $n_1 \ge n_0$ so that 
\[ |\nu_i |(A_m(n_1) \sm \{(\omega,0)\})<\eps \mbox{ for every } i<i_0 \text{ and } m < 3.\]
Applying Lemma \ref{main1} to $f=\chi_{A_0(n_1)}$ with $s\lesssim c$ (see Lemma \ref{p:2}); we obtain 
\[ t\ge 2\frac{c-\theta_0}{c-1}-\theta_0+\theta_1+\theta_2+a.\]
We also need an inequality suitable for the intermediate case when the values $\theta_i$'s are neither large nor small.
Suppose, for instance, that $1\notin J$, which means that $|\nu_i(A_1(n_1))| \lesssim c$.
Consider
\[h= \chi_{A_1(n_1)} - (1/2)\chi_{A_0(n_1)} - (1/2)\chi_{A_2(n_1)},\]
and note that 
\[\|h\| =\big\|(1/2) \chi_{A_1(n_1)} + (1/2)\big(\chi_{A_1(n_1)}-\chi_{A_0(n_1)}-\chi_{A_2(n_1)} \big) \big\| \lesssim c/2 + t/2,\]

Finally, the constraint of another type comes directly from Lemma \ref{sc:4}. For instance, if $J = \{0, 1, 2\}$, then
\[ \frac{t - c - 1}{c}(\theta_0+\theta_1+\theta_2)+a\ge 1.\]

Precise systems of equations are written in Appendix \ref{appendix}. 
We have approximated the minimal value of $t$ for which they have a solution in the following manner. 
First, we fix the value of $c$ to be a certain fraction of $t$ (such as $t/2$ or $(t+1)/2$).
If $c$ is linear with respect to $t$, these systems exhibit very notable property --- monotonicity: if there exists a solution for some $t_0$, then there is also one for any $t > t_0$.
This monotonicity allows us to run a simple binary search algorithm. 
From the known results, we have that the systems have no solutions for $t = 3$, while they do for $t = 5$.
Our binary search algorithm reveals that $c \approx (t+1/2)/2$ works best and our systems of equations have no solutions for $t \leq 3.53125$ allowing us to state the following theorem. 

\begin{theorem}
$\bmd\big(C([0,\omega]\times 3),C[0,\omega]\big) \geq 3.53125$.
\end{theorem}

The computations were performed in {\sc Mathematica}; the corresponding code is available on GitHub (see the file "Binsearch model.nb" in the \href{https://github.com/ememak/Bounds-for-Banach-Mazur-distance}{project files}).


\subsection{Upper bound}

Gergont and Piasecki in \cite{GP24}*{Section 3} introduced a very natural class of isomorphisms between $C(K)$ and $C(L)$ for $K=[0,\omega]\times 3$ and $L = [0,\omega]$. 
We briefly outline their construction here, mainly to provide an explicit value of the parameter $t$ that is optimal within this class.

It is interesting to note that just before this article was completed, the authors were informed that Marek Cuth from Prague had also obtained the same result.

\begin{theorem}
$\bmd\big(C([0,\omega]\times 3),C[0,\omega]\big) \leq \frac{4 + \sqrt[3]{73 - 6\sqrt{87}} + \sqrt[3]{73 - 6\sqrt{87}}}{3} \approx 3.87512...$
\end{theorem}

\begin{proof}
We define two matrices depending on a parameter $t$ with $3 \leq t \leq 4$, which will be the norm of the isomorphism we construct.
Put
{\large{\setstretch{1.2}
\[M = 
\begin{pmatrix}
t - 2           & -1                    & -1    \\
0               & t/2                   & -t/2   \\
\frac{t-2}{t}   & -\frac{t^2-5t+2}{4}   & -\frac{t^2-5t+2}{4}
\end{pmatrix},
\]}}
{\large{\setstretch{1.2}
\[C = 
\begin{pmatrix}
\frac{2t}{t+1} 	& 0						& 0\\
0				& \frac{t^2-t+2}{2t} 	& 0 \\
0 				& 0 					& \frac{t^2-t+2}{2t}
\end{pmatrix},
\]}}
denote the last row of $M$ by $M_3$ and let 
$M' = 
\begin{pmatrix}
M_3\\ M_3 \\ M_3
\end{pmatrix}$ be a 3x3 matrix with each row equal to $M_3$.
Given $f\in C(K)$, define
\[
\begin{pmatrix}
Tf(1)\\ Tf(2) \\ Tf(\omega)
\end{pmatrix} = 
M \cdot
\begin{pmatrix}
f(\omega, 0) \\ f(\omega, 1) \\ f(\omega, 2)
\end{pmatrix}
\]
and
\[
\begin{pmatrix}
Tf(3m) \\ Tf(3m+1) \\ Tf(3m+2)
\end{pmatrix} = 
C \cdot
\begin{pmatrix}
f(m, 0) \\ f(m, 1) \\ f(m, 2)
\end{pmatrix} + 
(M' - C) \cdot
\begin{pmatrix}
f(\omega, 0) \\ f(\omega, 1) \\ f(\omega, 2)
\end{pmatrix},
\]
where $m\in [1, \omega]$.
We use here a slightly modified notation compared with \cite{GP24}*{Section 3}, but $T$ still belongs to the same class of isomorphisms (in a slightly simplified form).

It is straightforward, though somewhat tedious (or best verified by computer), to check that

{\large{\setstretch{1.2}
\[M^{-1} = 
\begin{pmatrix}
\frac{t (t^2 - 5 t + 2)}{t^4 - 7 t^3 + 12 t^2 - 8 t + 8} & 0 & -\frac{4 t}{t^4 - 7 t^3 + 12 t^2 - 8 t + 8} \\
2/(t^3 - 5 t^2 + 2 t - 4) & 1/t & -(2 t)/(t^3 - 5 t^2 + 2 t - 4) \\
2/(t^3 - 5 t^2 + 2 t - 4) & -1/t & -(2 t)/(t^3 - 5 t^2 + 2 t - 4)
\end{pmatrix}
.\]}}
Next, define $S : C([1,\omega]) \to C([1,\omega] \times 3)$ as
\[
\begin{pmatrix}
Sg(\omega, 0) \\ Sg(\omega, 1) \\ Sg(\omega, 2)
\end{pmatrix} = 
M^{-1} \cdot
\begin{pmatrix}
g(1)\\ g(2) \\ g(\omega)
\end{pmatrix}
\]
and
\[
\begin{pmatrix}
Sg(m, 0) \\ Sg(m, 1) \\ Sg(m, 2)
\end{pmatrix} = 
C^{-1} \cdot
\left(
\begin{pmatrix}
g(3m)\\ g(3m+1) \\ g(3m+2)
\end{pmatrix} - 
(M' - C) \cdot M^{-1} \cdot
\begin{pmatrix}
g(1)\\ g(2) \\ g(\omega)
\end{pmatrix}
\right),
\]
for $g \in C([1,\omega])$ and $m\in \omega$. Then $S$ is an inverse of $T$.

Now the point is that if
\[t = \frac{4 + \sqrt[3]{73 - 6\sqrt{87}} + \sqrt[3]{73 - 6\sqrt{87}}}{3}\]
we have $\|T\| = t$ and $\|S\| = 1$.
\end{proof}
The analysis showing that the isomorphism constructed above is indeed optimal within the class defined by Gergont and Piasecki lies beyond the scope of this paper.

\appendix

\section{Four linear problems} \label{appendix}
In this appendix, we list the linear systems arising from the analysis outlined in Section~\ref{k3}. 
Recall that

\begin{enumerate}[(i)]
\item we can choose any $c\ge t/2$ (however, $c$ close to $t$ makes some cases trivial but others give weak bounds);
\item the problem is to determine the maximal value of $t$ for which none of these systems of inequalities has a solution
\end{enumerate}

Recall also that the code used to compute the corresponding values of $t$ was implemented in {\sc Mathematica} and is available on GitHub (see the file "Binsearch model.nb" in the \href{https://github.com/ememak/Bounds-for-Banach-Mazur-distance}{project files}).

\numberwithin{equation}{subsection}

\subsection{Case: \texorpdfstring{$J=\{0,1,2\}$.}{J = {0, 1, 2}.}}

\begin{equation}\label{7a}
t \ge \frac{2t}{t-1}+\theta_2 +a
\end{equation}
 
\begin{equation}\label{7b}
t \ge 2\frac{t-\theta_0}{t-1}-\theta_0+\theta_1+\theta_2 +a
\end{equation}

\begin{equation}\label{7c}
\frac{t - c - 1}{c}(\theta_0+\theta_1+\theta_2)+a \ge 1
\end{equation}
 
\begin{equation}\label{7d}
 0\le \theta_0\le \theta_1\le \theta_2, a\ge 0
\end{equation}

\subsection{Case: \texorpdfstring{$0\notin J$.}{0 not in J.}} 

\begin{equation}\label{6a}
 t\ge \frac{2t}{t-1} +\theta_2+a
\end{equation}
 
\begin{equation}\label{6b}
 t\ge 2\frac{c-\theta_0}{c-1}-\theta_0+\theta_1+\theta_2+a
\end{equation}
 
\begin{equation}\label{6b2}
 t\ge 2\frac{t/2+c/2-(\theta_0-(\theta_1+\theta_2)/2)}{t/2 + c/2 - 1}-\theta_0+\theta_1+\theta_2+a
\end{equation}

\begin{equation}\label{6c}
 \theta_0+\theta_1+\theta_2+ a\ge 1
\end{equation}
 
\begin{equation}\label{6d}
 0\le \theta_0\le \theta_1\le \theta_2, a\ge 0
\end{equation}

\subsection{Case: \texorpdfstring{$0\in J, 1\notin J$.}{0 in J, 1 not in J.}} 

\begin{equation}\label{8a}
 t\ge \frac{2t}{t-1} +\theta_2+a
\end{equation}
 
\begin{equation}\label{8b}
 t\ge 2\frac{t-\theta_0}{t-1}-\theta_0+\theta_1+\theta_2+a
\end{equation}
 
\begin{equation}\label{8c}
 t\ge 2\frac{c-\theta_1}{c-1}-\theta_1+\theta_0+\theta_2+a
\end{equation}
 
\begin{equation}\label{8d}
 t\ge 2\frac{t/2+c/2-(\theta_1-(\theta_0+\theta_2)/2)}{t/2+c/2-1}-\theta_1+\theta_1+\theta_2+a
\end{equation}
 
\begin{equation}\label{8e}
 \frac{t-c-1}{c}\theta_0+\theta_1+\theta_2+ a\ge 1
\end{equation}

\begin{equation}\label{8f}
 0\le \theta_0\le \theta_1\le \theta_2, a\ge 0
\end{equation}
\subsection{Case: \texorpdfstring{$0,1\in J, 2\notin J$.}{0, 1 in J, 2 not in J.}} 

\begin{equation}\label{9a}
 t\ge \frac{2t}{t-1} +\theta_2+a
\end{equation}
 
\begin{equation}\label{9b}
 t\ge 2\frac{t-\theta_0}{t-1}-\theta_0+\theta_1+\theta_2+a
\end{equation}
 
\begin{equation}\label{9c}
 t\ge 2\frac{c-\theta_2}{c - 1}-\theta_2+\theta_0+\theta_1+a
\end{equation}
 
\begin{equation}\label{9d}
 t\ge 2\frac{t/2+c/2-(\theta_2-(\theta_0+\theta_1)/2)}{t/2+c/2-1}-\theta_2+\theta_0+\theta_1+a
\end{equation}
 
\begin{equation}\label{9e}
 \frac{t-c-1}{c}\big(\theta_0+\theta_1\big)+\theta_2+ a\ge 1
\end{equation}

\begin{equation}\label{9f}
 0\le \theta_0\le \theta_1\le \theta_2, a\ge 0
\end{equation}

\section{Improving Lemma \texorpdfstring{\ref{main1}}{3.1}}

The following result is an effort to improve Lemma \ref{main1}. We use the same notation as in Section \ref{basic}.

\begin{lemma} \label{impr:1}
Suppose that $x\in K^{(1)}$, $C\sub K$ is a clopen set containing $x$ and $f\in C(K)$ satisfies
\begin{enumerate}[(i)]
 \item $t \ge \|Tf\|=\sup_i |\nu_i(f)|>|\nu(f)|$;
 \item $f(x)=1=\|f\|$,
 \item $f\chi_C = f$.
\end{enumerate}
Then either
\[\mbox{\rm (\ref{impr:1}.a)}\quad t\ge 2\frac{t - |\nu(f)| + |\nu|(K\sm C)}{t-1} - |\nu(f)| +|\nu|(K\sm C).\]
or
\[\mbox{\rm (\ref{impr:1}.b)}\quad t\ge 2\frac{t + |\nu(f)| - |\nu|(K\sm C)}{t-1} + |\nu(f)| + |\nu|(K\sm C).\]
\end{lemma}

\begin{proof}
We first fix a sequence of isolated points $x_n\in C$ converging to $x$. 
Write $e_n\in C(K)$ for the characteristic function of $\{x_n\}$.

Fix $\eps>0$.
Since $C$ is clopen, we can fix a norm-one function $f' \in C(K)$ supported outside of $C$ satisfying $\nu(f') \geq |\nu|(K\sm C) - \eps$. 
Without loss of generality, we also assume that $\nu(f) \geq 0$ (see Lemma \ref{main:2 Composition with an involution}).

Consider the functions
\[ g_n=\frac{1}{t+\eps} \cdot f +\left(1-\frac{1-\eps}{t+\eps}\right)\cdot e_n - \frac{1}{t+\eps} \cdot f'.\]
Note that $f(x_n) > 1 - \eps$ for large $n$ and then $\|g_n\|\ge g_n(x_n)\ge 1$.

As $\nu_i(f)\to \nu(f)$ and $\nu_i(f')\to \nu(f')$, there is $i_0$ such that for every $i\ge i_0$ we have $|\nu_i(f)-\nu(f)|, |\nu_i(f')-\nu(f')|< \eps$. 
Then fix $N$ such that for every $n\ge N$ and every $i<i_0$ we have $|\nu_i(e_n)|<\delta$,
where $\delta$ will be specify in a while.
Since $\|f - f'\| \leq 1$, we infer that for every $i<i_0$ and $n\ge N$,
\[ |\nu_i(g_n)| < \frac{t}{t+\eps}+\left(1-\frac{1-\eps}{t+\eps}\right)\cdot\delta < 1 \]
whenever $\delta$ is small enough.
In other words, we have checked that the initial measures cannot norm $g_n$ for large $n$: 

\begin{scclaim}
For every $n\ge N$ there is $i = i(n) \ge i_0$ such that $|\nu_i(g_n)|\ge \|g_n\|$.
\end{scclaim}

Note that a measure of finite variation may have only finitely many large atoms; hence the sequence $i(n)$ is unbounded. 

Now, let us analyse the signs of $\nu_{i(n)}(e_n)$, $\nu_{i(n)}(f)$ and $\nu_{i(n)}(f')$. 
We have assumed that $\nu_{i(n)}(f')\approx |\nu|(K\sm C) \geq 0$ and $\nu_{i(n)}(f) \geq 0$, so the only mystery is the sign of $\nu_{i(n)}(e_n)$.

If $\nu_{i(n)}(e_n) > 0$, then we have the following:
\[ 1 \le \nu_{i(n)}(g_n) \lesssim \frac{\nu(f)}{t} + \frac{t-1}{t} \cdot \nu_{i(n)}(e_n) - \frac{\nu(f')}{t}, \mbox{ so}\]
\[ \nu_{i(n)}(e_n) \gtrsim \frac{t - \nu(f) + |\nu|(K\sm C)}{(t-1)}.\]

We know that $f(x_n) \to 1$ when $n \to \infty$, so it follows that $\|2e_n - f\| \lesssim 1$ and $|\nu(C)| \gtrsim |2\nu(e_n) - \nu(f)|$. Thus
\[ |\nu|(K) = |\nu|(C) + |\nu|(K\sm C) \geq 2\frac{t - \nu(f) + |\nu|(K\sm C)}{(t-1)} - |\nu(f)| +|\nu|(K\sm C).\]

Otherwise, when $\nu_{i(n)}(e_n) < 0$, we have
\[ -1 \ge \nu_{i(n)}(g_n) \gtrsim \frac{\nu(f)}{t} + \frac{t-1}{t} \cdot \nu_{i(n)}(e_n) - \frac{\nu(f')}{t}, \mbox{ so}\]
\[ \nu_{i(n)}(e_n) \lesssim -\frac{t + \nu(f) - |\nu|(K\sm C)}{t-1}.\]
Using similar arguments as above, we obtain
\begin{align*}
|\nu|(K) &= |\nu|(C) + |\nu|(K\sm C) \geq 2\left|-\frac{t + \nu(f) - |\nu|(K\sm C)}{t-1} - \nu(f) \right| + |\nu|(K\sm C) = \\
&= 2\frac{t + |\nu(f)| - |\nu|(K\sm C)}{t-1} + |\nu(f)| + |\nu|(K\sm C).
\end{align*}

\end{proof}

Unfortunately, Lemma \ref{impr:1} is difficult to use, as it only gives a dichotomy, which greatly increases the number of cases necessary to consider.
However, it seems very useful. 

Consider $K = [0, \omega] \times 3$. 
Some numerical evidence suggests that it might improve the lower bound of $\bmd \big(C(K), C[0,\omega]\big)$ to approximately $3.6$ or more. 
It also seems rather difficult to construct an isomorphism $T : C(K) \to C([0, \omega])$ of small distortion where for $f = \chi_{[N, \omega] \times \{0\}}$ and $N$ large enough (as in Section \ref{k3}) the inequality \ref{impr:1}(b) holds.
It seems possible that the inequality \ref{impr:1}(a) always has to hold in this case, which, e.g., would mean that 
\[\lim_{k\in \omega} \bmd\big(C([0, \omega] \times k), C([0, \omega])\big) = 2 + \sqrt{5}.\]

\textbf{Declaration of generative AI and AI-assisted technologies in the manuscript preparation process}

During the preparation of this work the authors used Writefull and ChatGPT in order to correct the text grammatically and syntactically. After using these tools, the authors reviewed and edited the content as needed and take full responsibility for the content of the published article.

\bibliography{refs}

\end{document}